\documentclass[12pt]{amsart}
\usepackage[colorlinks=true,citecolor=magenta,linkcolor=magenta]{hyperref}
\usepackage{amssymb,verbatim,amscd,amsmath,graphicx}
\usepackage{enumerate}
\usepackage{graphicx}
\usepackage{caption}
\usepackage{subcaption}
\usepackage{pdfsync}
\usepackage{color,colortbl}
\usepackage{fullpage}
\usepackage{bbm}
\usepackage{comment}
\usepackage{tikz}
\usepackage{multirow}
\numberwithin{equation}{section}
\newtheorem{theorem}{Theorem}[section]

\newtheorem{proposition}[theorem]{Proposition}
\newtheorem{corollary}[theorem]{Corollary}

\theoremstyle{definition}

\newtheorem{conjecture}[theorem]{Conjecture}
\newtheorem{def-prop}[theorem]{Definition-Proposition}

\newtheorem{remark}[theorem]{Remark}
\newtheorem{example}[theorem]{Example}

\newtheorem*{acknowledgement}{Acknowledgements}

\DeclareMathOperator{\reg}{reg}

\DeclareMathOperator{\Ass}{Ass}

\newcommand{\CC}{{\mathbb C}}

\newcommand{\PP}{{\mathbb P}}

\newcommand{\NN}{{\mathbb N}}

\newcommand{\XX}{{\mathbb X}}

\def\mm{{\mathfrak m}}

\def\pp{{\frak p}}

\def\1{{\bf 1}}
\def\0{{\bf 0}}

\begin{document}
	
\title{The Initial Degree of Symbolic Powers of Ideals of Fermat Configuration of Points}

\author{Th\'ai Th\`anh Nguy$\tilde{\text{\^e}}$n}
\address{Tulane University \\ Department of Mathematics \\
	6823 St. Charles Ave. \\ New Orleans, LA 70118, USA}
\address{University of Education, Hue University, 34 Le Loi St., Hue, Viet Nam}
\email{tnguyen11@tulane.edu}
\urladdr{https://sites.google.com/view/thainguyenmath}


\begin{abstract}
Let $n \ge 2$ be an integer and consider the defining ideal of the Fermat configuration of points in $\PP^2$: 
$I_n=(x(y^n-z^n),y(z^n-x^n),z(x^n-y^n)) \subset R=\CC[x,y,z]$. In this paper, we compute explicitly the least degree of generators of its symbolic powers in all unknown cases. As direct applications, we easily verify Chudnovsky's Conjecture, Demailly's Conjecture and Harbourne-Huneke Containment problem as well as calculating explicitly the Waldschmidt constant and (asymptotic) resurgence number.
\end{abstract}

\maketitle


\section{Introduction} \label{sec.intro}

Let $n \ge 2$ be an integer and consider the \textit{Fermat ideal} 
$$I_n=(x(y^n-z^n),y(z^n-x^n),z(x^n-y^n)) \subset R=\CC[x,y,z].$$

This ideal corresponds to Fermat arrangement of lines (or Ceva arrangement in some literature) in $\PP^2$, more precisely, the variety of $I_n$ is a reduced set of $n^2+3$ points in $\PP^2$ \cite{HaSeFermat}, where $n^2$ of these points form the intersection locus of the pencil of curves spanned by $x^n- y^n$ and $x^n-z^n$, while the other $3$ are the coordinate points $[1 : 0 : 0], [0 : 1 : 0]$ and $[0 : 0 : 1]$. This set of points is said to be the Fermat configuration of points, justifying the terminology of Fermat ideal. Fermat ideal has attracted a lot of attention recently in commutative algebra research since it appeared as the first example of the non-containment between the third symbolic power and the second ordinary power of a defining ideal of a set of points in $\PP^2$, in the work of Dumnicki, Szemberg and Tutaj-Gasínska \cite{counterexamples} (when $n=3)$ and were generalized by Harbourne and Seceleanu \cite{HaSeFermat}. It is worth to emphasize that this is a quite surprising relation between ordinary and symbolic powers of ideal of points in $\PP^2$, since $I^2$ always contains $I^{(3)}$ where $I$ is the ideal of a general set of points \cite{BoH}. \par
\vspace{0.5em}

Since then, much has become known about Fermat ideals for $n\ge 3$. Fermat ideals can also be thought of the ideals determining the singular loci of the arrangements of lines given by the monomial groups $G(n,n,3)$, see \cite{DrabkinSeceleanu} or \cite{JustynaFermat}, \textit{Waldschmidt constant} and \textit{(asymptotic) resurgence number} of Fermat ideals have been computed in \cite{DHNSST2015} for $n\ge 3$. Nagel and Seceleanu in \cite{NagelSeceleanu} studied Rees algebra and symbolic Rees Algebra of Fermat ideals as well as the minimal generators, minimal free resolutions of all their ordinary powers and many symbolic powers. Specifically, it was shown that symbolic Rees algebra of $I_n$ is Noetherian; Castelnuovo-Mumford regularity of powers of $I_n$ and their reduction ideals were provided. In term of minimal generators, when $n\ge 3$, they provided the minimal generators for all ordinary powers as well all multiple of $n$ symbolic
powers. All known results are:
\begin{itemize}
    \item $\alpha(I_n^{(3k)})=3nk$ for all $k \ge 1$ by \cite[Theorem 2.1]{DHNSST2015}.
    \item $\alpha(I_n^{(nk)}) = n^2k$ for all $k\ge 1$ by \cite[Theorem 3.6]{NagelSeceleanu}.
    \item $\alpha(I_3^{(3m+2)})=9m+8$ for all $m\ge 0$ by \cite[Example 4.4]{MSS2018}.
\end{itemize}

In this paper, we will compute explicitly the least degree of generators (initial degree) of symbolic powers of $I_n$ for \textit{all} remaining cases as a contribution to provide a more complete picture of Fermat ideals. We summarize our results and results known in the literature in the following theorem.

\begin{theorem}
Let $n \ge 2$ be an integer and $I_n=(x(y^n-z^n),y(z^n-x^n),z(x^n-y^n))$ in $\CC[x,y,z]$ be a Fermat ideal. Then: 
\begin{enumerate}
    \item For $n\ge 2$, $\alpha(I_n^{(2)})=\alpha(I_n^{2})=2(n+1)$. (Theorem \ref{thm.secondpower} and Theorem \ref{thm.I2})
    \item $\alpha(I_2^{(2k)})=5k$ for all $k\ge 2$. (Theorem \ref{thm.I2})
    \item $\alpha(I_2^{(2k+1)})=5k+3$ for all $k\ge 0$. (Theorem \ref{thm.I2})
    \item $\alpha(I_3^{(3k)})=9k$ for all $k\ge 1$. (\cite[Theorem 2.1]{DHNSST2015})
    \item $\alpha(I_3^{(3k+1)})=9k+4$ for all $k\ge 1$. (Theorem \ref{thm.3m+1})
    \item $\alpha(I_3^{(3k+2)})=9k+8$ for all $k\ge 1$. (\cite[Example 4.4]{MSS2018})
    \item $\alpha(I_4^{(m)})=4m$ for all $m\ge 3$ and $m\not= 5$, $\alpha(I_4^{(5)})=21$ (Theorem \ref{thm.n=4})
    \item For $n \ge 5$, $\alpha(I_n^{(m)})=mn$ for all $m\ge 3$. (Theorem \ref{thm.n>5})
\end{enumerate}
\end{theorem}

\begin{table}[h!]
\centering
\tiny{
  \begin{tabular}{|p{1cm}||p{0.8cm}|p{0.8cm}|p{0.8cm}||p{0.8cm}|p{0.8cm}|p{0.8cm}|p{0.8cm}||p{0.8cm}|p{0.8cm}|p{1.2cm}||p{1.2cm}|p{0.8cm}|}
    \hline
    $n$ &
      \multicolumn{3}{c||}{2}&
      \multicolumn{4}{c||}{3}&
      \multicolumn{3}{c||}{4}&
      \multicolumn{2}{c|}{$n \ge 5$} \\
      \hline
    $m$ & 2 & $2k$ & $2k+1$ & 2 & $3k$ & $3k+1$ & $3k+2$ & 2 & 5 & $\ge 3$, $\not= 5$ & 2 & $\ge 3$ \\
    \hline
    $\alpha(I_n^{(m)})$ & 6 & $5k$ & $5k+3$ & 8 & $9k$ & $9k+4$ & $9k+8$ & 10 & 21 & $4m$ & $2(n+1)$ & $nm$ \\
    \hline
    $\widehat{\alpha}(I_n)$ &
      \multicolumn{3}{c||}{5/2}&
      \multicolumn{4}{c||}{3}&
      \multicolumn{3}{c||}{4}&
      \multicolumn{2}{c|}{$n$} \\
    \hline
    $\widehat{\rho}(I_n)$&
        \multicolumn{3}{c||}{$6/5$}&
      \multicolumn{9}{c|}{$(n+1)/n$} \\
    \hline
    $\rho(I_n)$ &
      \multicolumn{3}{c||}{$6/5$}&
      \multicolumn{9}{c|}{$3/2$} \\
    \hline
  \end{tabular}}
  \caption{The initial degrees and other invariants related to symbolic powers of $I_n$}
\end{table}

The above results give complete answer to the question of the least generating degree all symbolic powers of $I_n$. That includes the calculation for the ideal $I_2$, which is less considered in the aforementioned works. Note that the ideal $I_2$ is the ideal determining the singular locus of the arrangement of lines given by the pseudoreflection group $D_3$, see \cite{DrabkinSeceleanu}. It is worth to point out that the irregular value of $\alpha(I_4^{(5)})=21$. \par
\vspace{0.5em}

The \textit{containment problem} for an ideal $I$ is to determine the set of pairs $(m,r)$ for which $I^{(m)}\subseteq I^r$. The deep results in \cite{ELS, comparison, MaSchwede} show that $I^{(m)}\subseteq I^r$ whenever $m\ge Nr$. In order to characterize the pairs $(r, m)$ numerically, the \textit{resurgence} $\rho(I)$ is introduced in \cite{BoH}, and the \textit{asymptotic resurgence} $\widehat{\rho}(I)$ is introduced in \cite{AsymptoticResurgence}. It is known that for $n\ge 3$, $\rho(I_n)=\dfrac{3}{2}$, $\widehat{\rho}(I_n)=\dfrac{n+1}{n}$ and $\widehat{\alpha}(I_n)=n$ by \cite[Theorem 2.1]{DHNSST2015}, where $\widehat{\alpha}(I)$ is the \textit{Waldschmidt constant} of $I$. We will compute the Waldschmidt constant and (asymptotic) resurgence number of $I_2$. 

\begin{theorem}
For the ideal $I_2$, the Waldschmidt constant is $\widehat{\alpha}(I_2) = \dfrac{5}{2}$; the resurgence number and asymptotic resurgence number are $\rho(I_2)=\widehat{\rho}(I_2)=\dfrac{6}{5}$.
\end{theorem}

In an effort to improve the containment $I^{(m)}\subseteq I^r$ for $m\ge Nr$ as well as to deduce Chudnovsky's Conjecture, Harbourne and Huneke in \cite{HaHu} conjectured that the defining ideal $I$ for any set of points in $\PP^N$ satisfies some stronger containment, namely, $I^{(Nm)} \subseteq \mm^{m(N-1)}I^m$ and $I^{(Nm-N+1)} \subseteq \mm^{(m-1)(N-1)}I^m,$ for all $m \geqslant 1$. An interesting fact about the Fermat ideals is the verification and failure of the containment can be checked purely by their numerical invariants including the least degree of generators of symbolic powers, the regularity, or the resurgence number and the maximal degree of generators. As a consequence of the above computations, we will easily deduce that all Fermat ideals verify \textit{Harbourne-Huneke Containment}, \textit{stable Harbourne Containment} and some stronger containment. In \cite[Example 3.5]{BGHN2020}, we showed the stronger containment (which implies both Harbourne-Huneke containment) $I_n^{(2r-2)} \subseteq \mm^{r}I_n^r$ for $n\ge 3$ and for all $r \gg 0$. Here, we show the containment for all possible cases.

\begin{corollary}
For every $n\ge 2$, Fermat configuration ideals verifies the following containment:
\begin{enumerate}
    \item Harbourne-Huneke containment (see \cite[Conjecture 2.1]{HaHu}) $$I_n^{(2r)} \subseteq m^rI_n^r, \quad \forall r\ge 1.$$
    \item Harbourne-Huneke containment (see \cite[Conjecture 2.1]{HaHu}) $$I_n^{(2r-1)} \subseteq m^{r-1}I_n^r$$ $\forall r\ge 3$ if $n\ge 3$ and $\forall r \ge 1$ if $n=2$.
    \item A stronger containment $I_n^{(2r-2)} \subseteq m^{r}I_n^r, \quad \forall r\ge 5$.
\end{enumerate}
\end{corollary}

This work can also be thought of a modest contribution to the theory of Hermite interpolation. Specifically, given a set of points $Z = \{ P_1, \ldots, P_s\}$ in a projective space $\PP^n$ and positive integers $m_1,\ldots, m_s$, a fundamental problem in the theory of Hermite interpolation is to determine the least degree of a homogeneous polynomial that vanishes to order $m_i$ at the point $P_i$ for every $i=1,\ldots ,s$. This is a very hard problem even when $m_1=\ldots =m_s$. In this case when all $m_i$ are the same, thanks to the Zariski-Nagata theorem, the above interpolation problem is the same as asking for the least degree of a nonzero homogeneous polynomial in the $k-$th symbolic power of the defining ideal of $Z$ where $k=m_1=\ldots =m_s$. Hence, this work combined with many previous works provide a complete answer to the above question for $Z$ be a Fermat configuration of points in $\PP^2$.\par
\vspace{0.5em}
In order to understand the generating degrees of symbolic powers of $I_n$, we discuss the maximal degree of a set of minimal generators of symbolic powers of $I_n$, denoted by $\omega(I_n^{(m)})$. From the description of generating sets in \cite{NagelSeceleanu}, it can be seen that for $n\ge 3$, $\omega(I_n^{(m)}) = m(n+1)$, for $m=kn$ or $m=n-1$. By relating with another invariant $\beta(I_n^{(m)})$ that is defined in \cite[Definition 2.2]{HaHu}, we show that $\beta(J) \le \omega(J)$ for an ideal of points $J$ in general and use this to show that for $n\ge 3$, $\omega(I_n^{(m)})= m(n+1)$ for $m \ge n^2-3n+2$.\par
\vspace{0.5em}

It is worth to remark that in the works \cite{Malara2017FermattypeCO}, \cite{MalaraSzpond}, Malara and Szpond also studied the generalization of Fermat configuration in higher dimension to provide more counterexamples to the containment $I^{(3)} \subseteq I^2$. It turns out that these ideals share some similar properties to Fermat ideals. We will investigate these ideals in $3-$dimensional space with the same questions in the continuation paper \cite{ThaiFermatPlane} in order to keep this paper more concise and focused. We work over the field of complex numbers but our results hold over any algebraically closed field of characteristic $0$.

\begin{acknowledgement}
The author would like to thank his advisor, Tài Huy Hà, for introducing him this subject and giving many helpful suggestions and comments. He also thanks Alexandra Seceleanu and Ben Drabkin for comments on an early draft of this manuscript.	Finally, he thanks the referees for a careful reading of the paper and for many valuable suggestions. Results in the paper form a part of the author's thesis \cite{Thaithesis}.
\end{acknowledgement}

\section{Preliminaries}

Let $R = \CC[x_0,\ldots ,x_N]$ be the homogeneous coordinate ring of $\PP^N$, and let $\mm$ be its maximal homogeneous ideal. For a homogeneous ideal $I \subseteq R$, let $\alpha(I)$ denote the least degree of a nonzero homogeneous polynomial in $I$, and let 
$$I^{(m)} := \bigcap_{\pp \in \Ass(R/I)} I^m R_\pp \cap R$$
denote its $m$-th \emph{symbolic power}.\par
\vspace{0.5em}
Geometrically, given a set of distinct points $\XX \subseteq \PP^N$ and an integer $m \geqslant 1$, by the Zariski--Nagata Theorem \cite{Zariski,Nagata,EisenbudHochster} (cf. \cite[Proposition 2.14]{DDGHN}), the least degree of a nonzero homogeneous polynomial in the homogeneous coordinate ring $\CC[x_0,\ldots ,x_N]$ that vanishes at each point in $\XX$ of order at least $m$ is $\alpha ( I_\XX^{(m)} )$, where $I_\XX \subseteq \CC[x_0,\ldots ,x_N]$ is the defining ideal of $\XX$.\par
\vspace{0.5em}

The \textit{Waldschmidt constant} of $I$ is defined to be the limit and turned out to be the infimum 
$$\widehat{\alpha}(I):= \lim_{m\rightarrow \infty} \dfrac{\alpha(I^{(m)})}{m} = \inf _{m\rightarrow \infty} \dfrac{\alpha(I^{(m)})}{m}.$$ 
There is a tight connection between Waldschmidt constant and an algebraic manisfestion of the Seshadri constant, especially for a set of very general points, see \cite[Section 8]{Seshadri}.\par
\vspace{0.5em}

In studying the lower bound for the least degree of a homogeneous polynomial vanishing at a given set of points in $\PP^N$ with a prescribed order, Chudnovsky \cite{Chudnovsky1981} made the following conjecture.

\begin{conjecture}[Chudnovsky]
	\label{conj.Chud}
	Let $I$ be the defining ideal of a set of points $\XX\subseteq \PP^N_\CC$. Then, for all $n \geqslant 1$,
$$\dfrac{\alpha(I^{(n)})}{n} \geqslant \dfrac{\alpha(I)+N-1}{N}.$$ 
\end{conjecture}

Chudnovsky's Conjecture has been investigated extensively, for example, in \cite{EsnaultViehweg, BoH, HaHu, GHM2013, Dumnicki2015, DTG2017, FMX2018, BGHN2020}. Recently, the conjecture was proved for a \emph{very general} set of points \cite{DTG2017, FMX2018}, for a \emph{general} set of sufficiently many points \cite{BGHN2020} and recently, for any number of general points in \cite{SankhoThaiChudnovsky}. The conjecture was also generalized by Demailly \cite{Demailly1982} as follows.

\begin{conjecture}[Demailly] \label{conj.Demailly}
	Let $I$ be the defining ideal of a set of points $\XX \subseteq \PP^N_\CC$ and let $m \in \NN$ be any integer. Then, for all $n \geqslant 1$,
$$\dfrac{\alpha(I^{(n)})}{n} \geqslant \dfrac{\alpha(I^{(m)}) + N-1}{m+N-1}.$$
\end{conjecture}

Demailly's Conjecture for $N = 2$ was proved by Esnault and Viehweg \cite{EsnaultViehweg}. Recent work of Malara, Szemberg and Szpond \cite{MSS2018}, and of Chang and Jow \cite{CJ2020}, showed that for a fixed integer $m$, Demailly's Conjecture holds for a \emph{very general} set of sufficiently many points and for a \emph{general} set of $k^N$ points. In \cite{BGHN2-2020}, the results was extended for a \emph{general} set of sufficiently many points.\par
\vspace{0.5em}

The \textit{containment problem} for ideal $I$ is to determine the set $S_I$ of pairs $(r, m)$ for which $I^{(m)}\subseteq I^r$. The deep results in \cite{ELS, comparison, MaSchwede} show that $I^{(m)}\subseteq I^r$ whenever $m\ge Nr$, hence, $\{ (m,r) : m\ge Nr \} \subseteq S_I$. In order to characterize $S_I$ numerically, the \textit{resurgence number} $\rho(I)$ is introduced in \cite{BoH} as
$$\rho(I):= \sup \{ \dfrac{m}{r} : I^{(m)} \not\subseteq I^r \}$$
and the \textit{asymptotic resurgence number} $\widehat{\rho}(I)$ is introduced in \cite{AsymptoticResurgence} as
$$\widehat{\rho}(I):= \sup \{ \dfrac{m}{r} : I^{(mt)} \not\subseteq I^{rt} \text{  for  } t\gg 0 \}.$$
There are only few cases for which $S_I$ is known completely or the resurgence number have been determined. In general, $\rho(I)$ and $\widehat{\rho}(I)$ are different. \par
\vspace{0.5em}

In an effort to improve the containment $I^{(m)}\subseteq I^r$ for $m\ge Nr$, Harbourne and Huneke in \cite{HaHu} conjectured that the defining ideal $I$ for any set of points in $\PP^N$ satisfies some stronger containment, namely, $I^{(Nm)} \subseteq \mm^{m(N-1)}I^m$ and $I^{(Nm-N+1)} \subseteq \mm^{(m-1)(N-1)}I^m,$ for all $m \geqslant 1$. Knowing these containment clearly helps us to know the set $S_I$ and the numbers $\rho(I),\widehat{\rho}(I)$. Conversely, knowledge about $S_I,\rho(I),\widehat{\rho}(I)$ can be helpful to prove the containment. One useful result that we will use lately is the following: If $\dfrac{m}{r} > \rho(I)$ then by definition $I^{(m)} \subseteq I^r$ and suppose in addition, $\alpha(I^{(m)}) \ge a+ \omega(I^r)$ for some integer $a$, where $\omega(I)$ is the maximum degree of generators in a set of minimal generators of $I$, then $I^{(m)} \subseteq \mm^a I^r$. We refer interested readers to \cite{CHHVT2020} for more information about the Waldschmidt constant, resurgence number, containment between symbolic and ordinary powers of ideals.

\section{Fermat Ideals for $n\ge 3$}\label{sec.n>3}

In this section, we focus on the Fermat ideals $I_n=(x(y^n-z^n),y(z^n-x^n),z(x^n-y^n))$ for $n\ge 3$. Let us first recall some known results about degree of generators, Waldschmidt constants and (asymptotic) resurgence numbers of $I_n$. 

\begin{itemize}
    \item $\alpha(I_n^{(3k)})=3nk$ for all $k \ge 1$ by \cite[Theorem 2.1]{DHNSST2015}.
    \item $\alpha(I_n^{(nk)}) = n^2k$ for all $k\ge 1$ by \cite[Theorem 3.6]{NagelSeceleanu}.
    \item $\alpha(I_3^{(3m+2)})=9m+8$ for all $m\ge 0$ by \cite[Example 4.4]{MSS2018}.
    \item $\rho(I_n)=\dfrac{3}{2}$ and $\widehat{\rho}(I_n)=\dfrac{n+1}{n}$ by \cite[Theorem 2.1]{DHNSST2015}
\end{itemize}

Let $f_n=y^n-z^n,g_n=z^n-x^n,h_n=x^n-y^n$. Then $$I_n=(xf_n,yg_n,zh_n)=(f_n,g_n)\cap (x,y) \cap (y,z) \cap (z,x).$$
It is well-known that, since $f_n$ and $g_n$ form a regular sequence, for any $m\ge 1$ we have $$I_n^{(m)} = (f_n,g_n)^m\cap (x,y)^m \cap (y,z)^m \cap (z,x)^m.$$
Geometrically, recall that the Fermat configuration consists of $n^2+3$ points which are all points having each coordinate equal to a $n^{th}$ root of $1$ and the points $[0:0:1],[0:1:0],[1:0:0]$. Furthermore, these points are intersections of $3n$ lines $L_j$ which have equations: $$x-\epsilon^k y=0, y- \epsilon^k z=0, z- \epsilon^k x=0$$ for $k=0,1,\ldots , n-1$. Each of these lines contains exactly $n+1$ points $P_i$ (one coordinate point and $n$ other points), and each of the points $[0:0:1],[0:1:0],[1:0:0]$ is on exactly $n$ lines while each of the other $n^2$ points is on exactly $3$ lines.\par
\vspace{0.5em}
We will use these descriptions to explicitly compute the least degree of generators of $I_n^{(m)}$ in all remaining unknown cases, provided the knowledge about the Waldschmidt constant of $I_n$. This is our main strategy in this paper and the continuation paper \cite{ThaiFermatPlane}, we study a subsequence of $\alpha(I^{(m)})$, which gives us information about $\widehat{\alpha}(I)$, then use this to calculate other $\alpha(I^{(m)})$. \par
\vspace{0.5em}

The following $4$ theorems provide us all remaining unknown initial degrees of symbolic powers of $I_n$:

\begin{theorem}\label{thm.n>5}
For $n\ge 5$, we have $$\alpha(I_n^{(m)})=nm$$ for all $m\ge 3$.
\end{theorem}

\begin{proof}
First, for $3\le m \le n$, we observe that $$f_ng_nh_n(f_n,g_n)^{m-3} \subseteq I_n^{(m)} = (f_n,g_n)^m \cap (x,y)^m \cap (y,z)^m \cap (z,x)^m.$$
In fact, since $h_n=-(f_n+g_n)$, we have that $f_ng_nh_n(f_n,g_n)^{m-3} \subseteq (f_n,g_n)^m$. Also, it is clear that the sum of degree with respect to $x$ and degree with respect to $y$ of any monomials of $f_ng_nh_n$ is at least $n\ge m$ so $f_ng_nh_n \in (x,y)^m$. Similarly, $f_ng_nh_n\in (y,z)^m \cap (z,x)^m$.\par
\vspace{0.5em}
Thus, for $3\le m \le n$, $\alpha(I_n^{(m)}) \leq nm$. Since $\widehat{\alpha}(I_n)=n$, we have $\alpha(I_n^{(m)}) \geq nm$, therefore $\alpha(I_n^{(m)}) = nm$.\par
\vspace{0.5em}

Now for any $k\ge 2$, we claim that for $0\leq a \leq n-1$ $$(f_ng_nh_n)^k(f_n,g_n)^{k(n-3)-a} \subseteq I^{(kn-a)}.$$
The argument is identical to that of the above, the only notice here is that for $n\ge 5$, 
$$k(n-3)\ge 2(n-3) \ge n-1 \ge a.$$ 
Thus, for all $k\ge 2$ and $0\leq a \leq n-1$,
$$\alpha(I_n^{(kn-a)}) \leq (kn-a)n$$ hence, $\alpha(I_n^{(kn-a)}) = (kn-a)n$. Since for $m>n$, there are unique $k\ge 2$ and $0\le a \le n-1$ such that $m=kn-a$, we have that $\alpha(I_n^{(m)})=mn$ for all $m>n$. 
\end{proof}

\begin{theorem}\label{thm.n=4}
When $n=4$, $\alpha(I_4^{(5)})=21$. and for all $m\ge 3$ but $m \not = 5$, $$\alpha(I_4^{(m)})=4m.$$
\end{theorem}
\begin{proof}

With the same argument we have that:
\begin{enumerate}
    \item For $3\le m \le 4$, $f_4g_4h_4(f_4,g_4)^{m-3} \subseteq I_4^{(m)}$. The statement is true for $m=3,4$.
    \item For $k=2$ and $0\leq a \leq 2$, $$(f_4g_4h_4)^k(f_4,g_4)^{k(n-3)-a} \subseteq I_4^{(kn-a)}$$ i.e, 
    $$(f_4g_4h_4)^2(f_4,g_4)^{2-a}\subseteq I_4^{(8-a)}$$ for $0\le a \le 2$. Thus, the statement is true for $m=6,7,8$.
    \item For $k\ge 3$, we have $k(4-3)\ge 3$ so for $0\le a \le 3$, $$(f_4g_4h_4)^k(f_4,g_4)^{k(n-3)-a} \subseteq I_4^{(kn-a)}$$ therefore, the statement is true for all $m\ge 9$.
\end{enumerate}
In case $(2)$, the argument does not work for $m=5$ (since $a$ would be $3$). However, by the same argument, we can check that $zf_4^2g_4h_4^2 \in I_4^{(5)}$. Now suppose that $\alpha(I_4^{(5)}) \le 20$. Then there is a divisor $D$ of degree $20$ vanishing to order at least $5$ at every $P_i$ in the Fermat configuration. Since the intersection of $D$ and any line $L_j$ in the Fermat line arrangement consists of $5$ points to order at least $5$, by Bezout Theorem, each $L_j$ is a component of $D$ because $\deg(D).\deg(L_j)=20<5\cdot 5$. Moreover, in the Fermat configuration, each of the coordinate points is on exactly $4$ lines and each of the other points is on exactly $3$ lines. Hence, the divisor $D'=D- \sum_{j=1}^{12} L_j$ of degree $8$ vanishes to order at least $1$ at each coordinate point and to order at least $2$ along the others $16$ points. Now intersecting $D'$ with any of the lines $L_j$, again, since each lines $L_j$ contains exactly one coordinate point and $4$ other points, and $8<1\cdot 1 +4\cdot 2 $, we conclude by Bezout Theorem that each $L_j$ is a component of $D'$. This is a contradiction since the number of lines is $12$. Therefore, $\alpha(I_4^{(5)})=21$.
\end{proof}

\begin{remark}
From theorems \ref{thm.n>5} and \ref{thm.n=4}, when $n \ge 4$, $\alpha(I_n^{(m)})=nm$ for all $m\ge 3$, except for $n=4,m=5$. This can be predicted by means of Bezout Theorem. More precisely, suppose that $D$ is a divisor of degree $nm$ that vanishes to order at least $m$ along $n^2+3$ points of the Fermat configuration. By similar argument using Bezout Theorem, since $nm < m(n+1)$, each line $L_j$ is a component of $D$. Hence, the divisor $D_1=D- \sum_{j=1}^{3n} L_j$ is of degree $n(m-3)$ and vanishes to order at least $m-n$ (assuming $m > n$) at each coordinate point and to order at least $m- 3$ along the others $n^2$ points. If the number of lines is at most the degree of $D'$, that is $3n \le (m-3)n$, then Bezout Theorem would not yield contradiction. Note that in this case, by Bezout Theorem again, since $n(m-3) < 1\cdot (m-n) + n(m-3)$, each $L_j$ is again a component of $D_1$. Repeating this argument $k$ times whenever possible, the divisor $D_k=D- k \sum_{j=1}^{3n} L_j$ is of degree $n(m-3k)$ and vanishes to order at least $m-kn$ at each coordinate point and to order at least $m- 3k$ along the others $n^2$ points. Bezout Theorem would yield contradiction if $n(m-3k) < 1\cdot (m-kn) + n(m-3k)$ (hence, each $L_j$ is a component of $D_k$) and $n(m-3k)<3n$ (the degree of $D_k$ is less than the number of $L_j$). These two inequalities are equivalent to $m-kn \ge 1$ and $m-3k \le 2$, by combining them, we have $k(n-3)\le 1$, which only happens when $n=4$ and $k=1$ (thus, $m=5$), as we saw earlier.
\end{remark}

\begin{theorem}\label{thm.secondpower}
For all $n\ge 3$, $\alpha(I_n^{(2)})=2(n+1)$.
\end{theorem}

\begin{proof}
We know that $\alpha(I_n^{(2)}) \le \alpha(I_n^2) =2(n+1)$. Now suppose that $\alpha(I_n^{(2)}) \le 2n+1$. Then there is a divisor $D$ of degree $2n+1$ vanishing to order at least $2$ at every $P_i$ in the Fermat configuration. 
Since the intersection of $D$ and any line $L_j$ in the Fermat line arrangement consists of $n+1$ points to order at least $2$, by Bezout Theorem, each $L_j$ is a component of $D$ because $\deg(D).\deg(L_j)=2n+1<2(n+1)$. This is a contradiction since there are $3n$ lines and $3n>2n+1=\deg(D)$ when $n\ge 3$. Therefore, $\alpha(I_n^{(2)})=2(n+1)$ for all $n\geq 3$.
\end{proof}

It is known that
$\alpha(I_3^{(3m+2)})=9m+8$ for all $m\ge 0$ from \cite[Example 4.4]{MSS2018} and $\alpha(I_3^{(3m)})=9m$ for $m\ge 1$ from \cite[Theorem 2.1]{DHNSST2015} when $n=3$. Now we compute the remaining case $\alpha(I_3^{(3m+1)})$.

\begin{theorem}\label{thm.3m+1}
$\alpha(I_3^{(3m+1)}) = 9m+4$ for all $m\ge 0$.
\end{theorem}

\begin{proof}
Proceed by the argument in \cite[Example 4.4]{MSS2018}, suppose that there is $m\ge 1$ such that $\alpha(I_3^{(3m+1)})\le 9m+3$. Then there is a divisor $D$ of degree $9m+3$ vanishing to order at least $3m+1$ at every point of $12$ points $P_i$ in the Fermat configuration. Intersecting $D$ with any of the $9$ lines $L_j$, since each lines $L_j$ contains exactly $4$ points and $9m+3<4(3m+1)$, we conclude by Bezout Theorem that each $L_j$ is a component of $D$. Hence, there exists a divisor $D'=D- \sum_{j=1}^9 L_j$ of degree $9(m-1)+3$ vanishing to order at least $3(m-1)+1$ at every point of $P_i$. Repeating this argument $m$ times we get a contradiction with $\alpha(I_3)=4$. Thus $\alpha(I_3^{(3m+1)})\ge 9m+4$ for all $m$.\par
\vspace{0.5em}
On the other hand, by degree argument, $f_3^mg_3^mh_3^{m+1}z\in (f_3,g_3)^m \cap (x,y)^m \cap (y,z)^m \cap (z,x)^m = I_3^{(3m+1)}$ so $\alpha(I_3^{(3m+1)})\le 9m+4$ for all $m$. Therefore, $\alpha(I_3^{(3m+1)})= 9m+4$ for all $m$.
\end{proof}

 \begin{example}
It is worth to point out that the first immediate application of the above calculations combining with the already known cases is the verification of Chudnovsky's Conjecture and Demailly's Conjecture, although the general case is already known from \cite{EsnaultViehweg}. For any $n\ge 3$, Fermat ideals verify
\begin{enumerate}
    \item Chudnovsky's Conjecture $\widehat{\alpha}(I_n)\ge \dfrac{\alpha(I_n)+1}{2}$.
    \item Demailly's Conjecture $\widehat{\alpha}(I_n)\ge \dfrac{\alpha(I_n^{(m)})+1}{m+1}$ for all $m\ge 1$.
\end{enumerate}
\end{example}

\begin{proof}
Directly from the formulae of $\widehat{\alpha}(I_n)$ and $\alpha(I_n^{(m)})$.
\end{proof}

The following containment are also direct consequences of the above calculations about $\alpha(I_n^{(m)})$. Note that these containment (and in fact, the stable version of them, i.e, the containment for $r \gg 0$ imply Chudnovsky's Conjecture). First, in \cite[Example 3.5]{BGHN2020}, we showed the stronger containment (which implies both Harbourne-Huneke containment) $$I_n^{(2r-2)} \subseteq \mm^{r}I_n^r$$ for $r=6$ and thus for all $r \gg 0$. In particular, from the proof of \cite[Theorem 3.1]{BGHN2020}, the containment hold for $r\ge 12^2=144$. Here we show that the containment hold for all $r\ge 5$. Notice that for $r\le 4$, since the resurgence number $\rho(I_n)=\dfrac{3}{2}$, we know that $I_n^{(2r-2)} \not \subseteq I_n^r$.

\begin{corollary}
For every $n\ge 3$, Fermat configuration ideal verifies the following containment $$I_n^{(2r-2)} \subseteq \mm^{r}I_n^r, \quad \forall r\ge 5.$$
\end{corollary}

\begin{proof}
As before, since for all $r\ge 5$,   $I_n^{(2r-2)} \subseteq I_n^r$, it suffices to check the inequalities $$\alpha(I_n^{(2r-2)})\ge r+\omega(I_n^r)$$
case by case.
\begin{enumerate}
    \item For $n=3$, we have $r+\omega(I_n^r)=5r$ and $\alpha(I_n^{(2r-2)})=
    \begin{cases}
    9m, & \text{if $2r-2=3m$}.\\
    9m+4, & \text{if $2r-2=3m+1$}\\
    9m+8, & \text{if $2r-2=3m+2$}.
  \end{cases}$
  \begin{enumerate}
      \item If $2r-2=3m$, the inequality becomes $9m\ge \dfrac{5}{2}(3m+2)$ which is equivalent to $3m\ge 10$. Since $r\ge 5$, we have $3m \ge 8$. Moreover, $3m=2r-2$ can't be $8$ or $9$. 
      \item If $2r-2=3m+1$, the inequality becomes $9m+4\ge \dfrac{5}{2}(3m+3)$ which is equivalent to $3m\ge 7$ (which is true because $3m=2r-3 \ge 7$).
      \item If $2r-2=3m+2$, the inequality becomes $9m+8\ge \dfrac{5}{2}(3m+4)$ which is equivalent to $3m\ge 4$.
  \end{enumerate}
  \item For $n\ge 4$, since $r\ge 5$, $2r-2 \ge 8$ so we have $\alpha(I_n^{(2r-2)})=(2r-2)n$, and the inequalities $(2r-2)n\ge r(n+1)+r \Leftrightarrow (r-2)n \ge 2r $, which is true since we have $(r-2)n\ge 4(r-2)\ge 2r$ for all $r\ge 5$.
\end{enumerate}
\end{proof}

Although the above containment imply the Harbourne-Huneke containment for $r\ge 5$, we can check easily that Harbourne-Huneke containment hold for all possible $r$ by our computations.

\begin{corollary}
For every $n\ge 3$, Fermat configuration ideal verifies Harbourne-Huneke containment (see \cite[Conjecture 2.1]{HaHu}) $$I_n^{(2r)} \subseteq \mm^rI_n^r, \quad \forall r \ge 1.$$
\end{corollary}

\begin{proof}
Since $I_n^{(2r)} \subseteq I_n^r, \forall r\ge 1$ the above containment come from the fact that $$\alpha(I_n^{(2r)})\ge r+\omega(I_n^r)$$ for all $n\ge 3$ and $r \ge 1$ (the case $r=0$ is trivial). Indeed, we check case by case
\begin{enumerate}
    \item For $n\ge 3$ and $r=1$, we have $2(n+1)\ge n+1+1$
    \item For $n=3$, we have $r+\omega(I_n^r)=5r$ and $\alpha(I_n^{(2r)})=
    \begin{cases}
    9m, & \text{if $2r=3m$}\\
    9m+4, & \text{if $2r=3m+1$}\\
    9m+8, & \text{if $2r=3m+2$}.
  \end{cases}$
    \item For $n\ge 4$ and $r\ge 2$, we have $r+\omega(I_n^r)=r(n+2)$ and $\alpha(I_n^{(2r)})= 2rn$. 
\end{enumerate}
\end{proof}

\begin{corollary}
For every $n\ge 3$, Fermat configuration ideal verifies Harbourne-Huneke containment (see \cite[Conjecture 4.1.5]{HaHu}) $$I_n^{(2r-1)} \subseteq \mm^{r-1}I_n^r, \quad \forall r\ge 3.$$
\end{corollary}

\begin{proof}
Since $\rho(I_n)=\dfrac{3}{2}$, for all $r\ge 3$,   $I_n^{(2r-1)} \subseteq I_n^r$, the above containment comes from the fact that $$\alpha(I_n^{(2r-1)})\ge r-1+\omega(I_n^r)$$ for all $n\ge 3$ and $r\ge 3$. Notice that when $r=1,r=2$, containment $I_n^{(2r-1)} \subseteq I_n^r$ fail. We check case by case
\begin{enumerate}
    \item For $n=3$, we have $r-1+\omega(I_n^r)=5r-1$ and $\alpha(I_n^{(2r-1)})=
    \begin{cases}
    9m, & \text{if $2r-1=3m$}.\\
    9m+4, & \text{if $2r-1=3m+1$}\\
    9m+8, & \text{if $2r-1=3m+2$}.
  \end{cases}$
  \item For $n=4,r=3$, we have $\alpha(I_n^{(5)})=21>3+5\cdot 3=r+\omega(I_n^3)$.
  \item For $n\ge 4, r\not= 3$, we have $\alpha(I_n^{(2r-1)})=(2r-1)n$, and the inequalities $(2r-1)n\ge r(n+1)+r-1 \Leftrightarrow (r-1)n \ge 2r-1 $, which is true since $(r-1)n\ge 4(r-1)\ge 2r-1$ for all $r\ge 3$.
\end{enumerate}
\end{proof}

We end this section by calculating another invariant related to generating degree of symbolic powers of $I_n$. As introduced in \cite[Definition 2.2]{HaHu}, for a homogeneous ideal $J \subset \CC[\PP^N]$, define $\beta(J)$ to be the smallest integer $t$ such that $[J]_t$ contains a regular sequence of length two where $[J]_t$ is the graded component of degree $t$ of $J$. It turns out that when $N=2$ and $J$ is a defining ideal of a finite set of (fat) points, $\beta(J)$ is in fact the least degree $t$ such that the zero locus of $[J]_t$ is $0$-dimensional since the condition that $[J]_t$ contains a regular sequence of length two is equivalent to the condition that all elements of $[J]_t$ does not have a non-constant common factor, c.f \cite{HaSeFermat}. This invariant is related to the maximum degree of generators in a set of minimal generators of $J$ as follows:

\begin{proposition}
Let $J \subset \CC[x,y,z]$ is a defining ideal of a set of finite (fat) points. Then $\beta(J) \le \omega(J)$.
\end{proposition}
\begin{proof}
Suppose that $J=\langle g_1, \ldots ,g_k \rangle$ and  $\omega(J) < \beta(J)$, then the zero locus of $[J]_{\omega(J)}$ is not $0-$dimensional by definition of $\beta(J)$. For $j=1,\ldots ,k$, consider the set $A$ consists of all forms $g_jx^{d_j},g_jy^{d_j},g_jz^{d_j}$ where $d_j =\omega(J)-\deg(g_j)$. Then $A\subset [J]_{\omega(J)}$, hence the zero locus of $A$ is not $0-$dimensional. This is a contradiction since the zero locus of $A$ is also the zero locus of $J$.
\end{proof}

In the following result, we calculate $\beta(I_n^{(m)})$ and get an immediate bound for $\omega(I_n^{(m)})$.

\begin{proposition}
For $n\ge 3$ and $m\ge 1$, $\beta(I_n^{(m)})=m(n+1)$.
\end{proposition}
\begin{proof}
First, for any $n\ge 3$, note that $I_n^m$ is generated by all generators of the same degree $m(n+1)$, hence, $\beta(I_n^m)=m(n+1)$ as $[I_n^m]_{m(n+1)}$ contains a regular sequence of length two. Since $I_n^m \subseteq I_n^{(m)}$, we have that $\beta(I_n^{(m)}) \le m(n+1)$ for all $m \ge 1$. On the other hand, for any $m\ge 1$, let $f \in [I_n^{(m)}]_t$ be any element where $t<m(n+1)$. Recall that each line $L_j$ in the Fermat configuration contains exactly $n+1$ points of the configuration. Intersecting any line $L_j$ with the variety defined by $f$, since $f$ vanishes at every point in the configuration to order at least $m$, by Bezout's theorem, since $\deg(f)\deg(L_j)<m(n+1)$, $L_j$ is a component of the variety of $f$. Therefore, for any $t<m(n+1)$, $L_j$ is a component of the zero locus of $[I_n^{(m)}]_t$, i.e, $\beta(I_n^{(m)}) \ge m(n+1)$. This ends the proof.
\end{proof}

\begin{corollary}
For $n\ge 3$ and $m\ge 1$, $\omega(I_n^{(m)})\ge m(n+1)$. Moreover, for each $n \ge 3$, $\omega(I_n^{(m)})= m(n+1)$ for all $m\ge n^2-3n+2$.
\end{corollary}

\begin{proof}
$\omega(I_n^{(m)})\ge m(n+1)$ for all $m$ is immediate. On the other hand, by  \cite[Theorem 3.10, Remark 3.11]{NagelSeceleanu}, $\reg(I_n^{(m)})=m(n+1)$ for all $n\ge 3$ and $m\ge n^2-3n+2$. Thus $\omega(I_n^{(m)})\le \reg(I_n^{(m)}) = m(n+1)$ for all $n\ge 3$ and $m\ge n^2-3n+2$.
\end{proof}

\begin{remark}
It is known that for any $k$ and $n\ge 3$, $\omega(I_n^{(kn)}) = \beta(I_n^{(kn)}) = kn(n+1)$, since any generators of $I_n^{(kn)}$ with degree less than $kn(n+1)$ must be divisible by $fgh$ and hence, none of two elements in degree less than $kn(n+1)$ of it form a regular sequence (because they always share the common factor $fgh$), see \cite[Theorem 3.7]{NagelSeceleanu}. In the same paper, it is also known that $\omega(I_n^{(n-1)})=\beta(I_n^{(n-1)})= (n-1)(n+1)$ with the same reason. As in above corrolary, $\omega(I_n^{(m)})=\beta(I_n^{(m)})= m(n+1)$ for $m\gg 0$. It is reasonable to ask if this is the case for all $m$. It is suggested by Macaulay2 \cite{M2} that in fact $\omega(I_n^{(m)})= m(n+1)$ for all $m\ge 1$.
\end{remark}

\section{Fermat Ideal $I_2$}
In this section, we will deal with the ideal $I_2=(x(y^2-z^2),y(z^2-x^2),z(x^2-y^2))$. Unlike the ideals $I_n$ for $n\ge 3$, this ideal satisfies the Harbourne containment $I_2^{(3)} \subseteq I_2^2$. This is probably one reason why it is less considered in the literature. We will see later that this ideal is very different to Fermat ideals when $n\ge 3$ in terms of generating degree of symbolic powers and hence, in terms of Waldschmidt constant and (asymptotic) resurgence number.\par
\vspace{0.5em}

In the following, we will compute the least degree of generators of symbolic powers of $I_2$ as well as its Waldschmidt constant, (asymptotic) resurgence number in order to complete the picture of Fermat ideals. This ideal $I_2$ is also known to be the ideal determining the singular locus of the arrangement of lines given by the pseudoreflection group $D_3$. In general, suppose that $G \subseteq \text{GL}_{n+1}(\CC)$ is a finite group generated by pseudoreflections, where a pseudoreflection is a non-identity linear transformation that fixes a hyperplane pointwise and has finite order. Geometrically, we can view the generators of $G$ as a hyperplane arrangement where the hyperplanes are pointwise fixed by the pseudoreflections of $G$. It is shown in \cite[Proposition 3.9]{DrabkinSeceleanu} that the singular locus (the Jacobian ideal) of the arrangement of lines correspond to $G(2,2,3)=D_3$ is given by $I_2$. \par
\vspace{0.5em}

First, recall that, geometrically, $I_2$ is the defining ideal of the singular locus of the line arrangement in $\PP^2$ that consists of $6$ lines $L_j$ whose equations are $$x=\pm y, y=\pm z, z=\pm x$$
These $6$ lines intersect at $7$ points $P_i$ which are $[1:0:0], [0:1:0], [0:0:1], [1:1:-1],$ $[1:-1:1], [1:-1:-1]$ and $[1:1:1]$ such that the first $3$ points lie on $2$ lines each and the rest lie on $3$ lines each; and each line contains exactly $3$ points.\par
\vspace{0.5em}

Similar to $I_n$ for $n\ge 3$, we can write 
$$ I_2= (x^2-y^2,y^2-z^2) \cap (x,y) \cap (y,z) \cap (z,x)$$ so that $$I_2^{(m)}=(x^2-y^2,y^2-z^2)^m \cap (x,y)^m \cap (y,z)^m \cap (z,x)^m \ \ \forall m\ge 1. $$

\begin{theorem}\label{thm.I2}
For ideal $I_2$ we have the following
\begin{enumerate}
    \item $\widehat{\alpha}(I_2) = \dfrac{5}{2}.$
    \item $\alpha(I_2^{(2k)})=5k$, for all $k\ge 2$.
    \item $\alpha(I_2^{(2k+1)})=5k+3$, for all $k\ge 0$.
    \item $\alpha(I_2^{(2)})=6.$
\end{enumerate}
\end{theorem}

\begin{proof}
By \cite[Theorem 2.3]{FaGHLMaS}, we can check that $\widehat{\alpha}(I_2) \ge \dfrac{5}{2}$. 
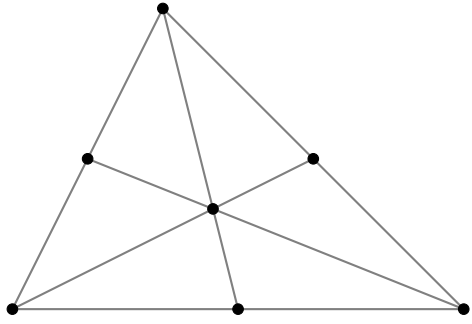
\begin{center}
\begin{tikzpicture}
\draw[gray, thick] (0,4) -- (-2,0);
\draw[gray, thick] (-2,0) -- (4,0);
\draw[gray, thick] (0,4) -- (4,0);
\draw[gray, thick] (0,4) -- (1,0);
\draw[gray, thick] (4,0) -- (-1,2);
\draw[gray, thick] (-2,0) -- (2,2);
\filldraw[black] (0,4) circle (2pt);
\filldraw[black] (-2,0) circle (2pt);
\filldraw[black] (4,0) circle (2pt);
\filldraw[black] (1,0) circle (2pt);
\filldraw[black] (-1,2) circle (2pt);
\filldraw[black] (2,2) circle (2pt);
\filldraw[black] (2/3,4/3) circle (2pt);
\end{tikzpicture}
\end{center}
In particular, we have that $\alpha(I_2^{(2k)}) \ge 5k$ and $\alpha(I_2^{(2k+1)}) \ge 5k+3$, for all $k\ge 1$. We will show the reverse by showing there exists some element with the desired degree in the symbolic powers.\par
\vspace{0.5em}
Indeed, we consider case by case. Denote $K=(x^2-y^2,y^2-z^2)$, we have \par
\vspace{0.5em}

\begin{enumerate}
    \item Case $1$: when $m=4k$. Consider $F=(x^2-y^2)^{2k}(y^2-z^2)^k(z^2-x^2)^kz^{2k}$ that has degree $10k$. Since $(x^2-y^2)^{2k}(y^2-z^2)^k(z^2-x^2)^k \in K^{4k}$, $(x^2-y^2)^{2k} \in (x,y)^{4k}$, $z^{2k}(y^2-z^2)^{k} \in (y,z)^{4k}$ and $z^{2k}(z^2-x^2)^{k} \in (z,x)^{4k}$ it follows that $$F\in K^{4k} \cap (x,y)^{4k} \cap (y,z)^{4k} \cap (z,x)^{4k}=I_2^{(4k)}, \ \forall k\ge 1.$$
    By similar argument we see that
    \item Case $2$: when $m=4k+2$. $F=(x^2-y^2)^{2k}(y^2-z^2)^{k+1}(z^2-x^2)^{k+1}xyz^{2k-1}$ has degree $10k+5$ and $$F\in K^{4k+2} \cap (x,y)^{4k+2} \cap (y,z)^{4k+2} \cap (z,x)^{4k+2}=I_2^{(4k+2)}, \ \forall k\ge 1.$$
    \item Case $3$: when $m=4k+1$. $F=(x^2-y^2)^{2k+1}(y^2-z^2)^{k}(z^2-x^2)^{k}z^{2k+1}$ has degree $10k+3$ and $$F\in K^{4k+1} \cap (x,y)^{4k+1} \cap (y,z)^{4k+1} \cap (z,x)^{4k+1}=I_2^{(4k+1)}, \ \forall k\ge 0.$$
    \item Case $4$: when $m=4k+3$. $F=(x^2-y^2)^{2k+1}(y^2-z^2)^{k+1}(z^2-x^2)^{k+1}xyz^{2k}$ has degree $10k+8$ and $$F\in K^{4k+3} \cap (x,y)^{4k+3} \cap (y,z)^{4k+3} \cap (z,x)^{4k+3}=I_2^{(4k+3)}, \ \forall k\ge 0.$$
\end{enumerate}

Thus, $\alpha(I_2^{(2k)})\le 5k$, for all $k\ge 2$ and $\alpha(I_2^{(2k+1)}) \le 5k+3$, for all $k\ge 0$. It follows that statements $(2)$ and $(3)$ are true, and since the Waldschmidt constant is the infimum of the initial degrees, $\widehat{\alpha}(I_2) \le \dfrac{\alpha(I_2^{(2k)})}{2k}= \dfrac{5}{2}$, hence, $(1)$ follows as well. Part $(4)$ can be checked directly by Macaulay2 or by Bezout theorem argument as follows: We know that $\alpha(I_2^{(2)}) \le \alpha(I_2^2) =6$. Now suppose that $\alpha(I_2^{(2)}) \le 5$. Then there is a divisor $D$ of degree $5$ vanishing to order at least $2$ at every point $P_i$. Since the intersection of $D$ and any line $L_j$ consists of $3$ points to order at least $2$, we get a contradiction to Bezout theorem because $\deg(D)\cdot \deg(L_j)=5<2\cdot 3$.
\end{proof}

For $I_2$, the asymptotic resurgence number and the resurgence number turn out to be the same.

\begin{theorem}
The resurgence number and asymptotic resurgence number of $I_2$ are: $$\rho(I_2)=\widehat{\rho}(I_2)=\dfrac{6}{5}.$$
\end{theorem}

\begin{proof}
The asymptotic resurgence number $\widehat{\rho}(I_2)=\dfrac{6}{5}$ follows from \cite[Theorem 1.2]{AsymptoticResurgence} that $$\dfrac{6}{5}=\dfrac{\alpha(I_2)}{\widehat{\alpha}(I_2)}\le \widehat{\rho}(I_2)\le \dfrac{\omega(I_2)}{\widehat{\alpha}(I_2)}=\dfrac{6}{5}.$$

By \cite[Theorem 2.5]{NagelSeceleanu}, since $I_2$ is a strict almost complete intersection ideal with minimal generators of degree $3$ and its module sygyzy is generated in degree $1$ and $2$, the minimal free resolution of $I_2^r$ is:
$$ 0 \to R(-3(r+1))^{{r} \choose 2}\stackrel{\psi}\to \begin{matrix} R(-3r-1)^{{r+1} \choose 2} \\
\oplus \\ R(-3r-2)^{{r+1} \choose 2}\end{matrix} \stackrel{\varphi}\to R(-3r)^{{r+2} \choose 2} \to I_2^r\to 0 $$
for any $r\ge 2$, in particular, $\reg(I_2^r)=3r+1$ for all $r\ge 2$, where $\reg(I)$ denotes the Castelnuovo-Mumford regularity of $I$. We also have $\reg(I_2)=4$. By \cite[Theorem 1.2]{AsymptoticResurgence} again, we have $\rho(I_2) \ge \dfrac{\alpha(I_2)}{\widehat{\alpha}(I_2)} =\dfrac{6}{5}$. Conversely, for any $\dfrac{m}{r}>\dfrac{6}{5}$, we have: 
\begin{itemize}
    \item If $m=2k$ then $10k>6r$ implies that $10k\ge 6r+2$ or $5k\ge 3r+1$ since both numbers are even. It follows that $\alpha(I_2^{(2k)}) \ge 5k \ge 3r+1 =\reg(I_2^r)$ and hence, $I_2^{(2k)} \subseteq I_2^r$.
    \item If $m=2k+1$ then $10k+2>6r$ implies that $10k+6\ge 6r+2$ or $5k+3\ge 3r+1$. It follows that $\alpha(I_2^{(2k+1)}) =5k+3 \ge 3r+1 =\reg(I_2^r)$ and hence, $I_2^{(2k+1)} \subseteq I_2^r$.
\end{itemize}
Thus, for any $\dfrac{m}{r}>\dfrac{6}{5}$, we have $I_2^{(m)} \subseteq I_2^r$, i.e, $\rho(I_2)=\dfrac{6}{5}$.
\end{proof}

\begin{example}
It is worth to point out that the first direct consequence of the above calculation is the verification of $I_2$ to Chudnovsky's Conjecture and Demailly's Conjecture, although the general case is already known from \cite{EsnaultViehweg}. Ideal $I_2$ verifies
\begin{enumerate}
    \item Chudnovsky's Conjecture $\widehat{\alpha}(I_2)\ge \dfrac{\alpha(I_2)+1}{2}$.
    \item Demailly's Conjecture $\widehat{\alpha}(I_2)\ge \dfrac{\alpha(I_2^{(m)})+1}{m+1}$ for all $m\ge 1$.
\end{enumerate}
\end{example}

\begin{proof}
Directly from the formulae of $\widehat{\alpha}(I_2)$ and $\alpha(I_2^{(m)})$.
\end{proof}

\begin{example}
Another difference between $I_2$ and $I_n$ when $n\ge 3$ is that while $I_n^{(nk)}=(I_n^{(n)})^k$ for all $k$ and $n\ge 3$ by \cite{NagelSeceleanu}, $I_2^{(4)} \not = (I_2^{(2)})^2$. In fact, as in the proof of theorem \ref{thm.I2}, $$F=(x^2-y^2)^2(y^2-z^2)(z^2-x^2)z^2 \in I_2^{(4)}$$ 
whereas, we can check that $F \not \in (I_2^{(2)})^2$. Moreover, we can also check that $I_2^{(6)} \not = (I_2^{(3)})^2$ since $$G=(x^2-y^2)^2(y^2-z^2)^2(z^2-x^2)^2xyz \in I_2^{(6)} \setminus (I_2^{(3)})^2.$$ 
It is suggested by Macaulay2 \cite{M2} that $I_2^{(8)} = (I_2^{(4)})^2$. It is interesting to know if $I_2^{(4k)}=(I_2^{(4)})^k$ for all $k$.
\end{example}

As all other Fermat ideals, $I_2$ also satisfies the following containment. In \cite[Example 3.7]{BGHN2020}, we showed the stronger containment (which implies both Harbourne-Huneke) containment $$I_2^{(2r-2)} \subseteq \mm^{r}I_2^r$$ for $r=5$ (by Macaulay2) and thus for all $r \gg 0$ by our method. In particular, from the proof of \cite[Theorem 3.1]{BGHN2020}, the containment hold for $r\ge 10^2=100$. Here we show that the containment hold for all $r\ge 5$. 

\begin{corollary}
For every $n\ge 3$, ideal $I_2$ verifies the following stronger containment $$I_2^{(2r-2)} \subseteq \mm^{r}I_2^r, \quad \forall r\ge 5.$$
\end{corollary}

\begin{proof}
As before, since $\rho(I_2)=\dfrac{6}{5}$, we know that $I_2^{(2r-2)} \subseteq I_2^r$ for $r\ge 3$. Hence, the containment follows from the inequalities $$\alpha(I_2^{(2r-2)})=5r-5 \ge r+3r = r+\omega(I_2^r)$$ for all $r\ge 5$.
\end{proof}

We can also detect the failure of the containment in the remaining cases by only using formulae for $\alpha(I_2^{(m)})$ as follows.

\begin{remark}
For $r\le 4$, the above containment fail. In fact, notice that for $r\le 2$, since $\rho(I_2)=\dfrac{6}{5}$, we know that $I_2^{(2r-2)} \not \subseteq I_2^r$. When $r=3$, since $\alpha(I_2^{(4)})=10 < 3+9 = 3+\alpha(I_2^3)$, we see that the containment $I_2^{(4)} \subseteq \mm^3I_2^3$ fails. Similarly, for $r=4$, since $\alpha(I_2^{(6)})=15 < 4+12 = 4+\alpha(I_2^4)$, we see that the containment $I_2^{(6)} \subseteq \mm^4I_2^4$ also fails.
\end{remark}

Although the above containment imply the Harbourne-Huneke containment for $r\ge 5$, we can check easily that Harbourne-Huneke containment hold for all possible $r$ by our computations.

\begin{corollary}
Ideal $I_2$ verifies Harbourne-Huneke containment (see \cite[Conjecture 2.1]{HaHu}) $$I_2^{(2r)} \subseteq \mm^rI_2^r, \quad \forall r \ge 1.$$
\end{corollary}

\begin{proof}
Since for all $r$, $I_2^{(2r)} \subseteq I_2^r$, the containment follow from the fact that $$\alpha(I_2^{(2r)}) \ge 5r\ge r+3r = r+\omega(I_2^r)$$
for all $r$.
\end{proof}

\begin{corollary}
Ideal $I_2$ verifies Harbourne-Huneke containment (see \cite[Conjecture 4.1.5]{HaHu}) $$I_2^{(2r-1)} \subseteq \mm^{r-1}I_2^r, \quad \forall r\ge 1.$$
\end{corollary}

\begin{proof}
Since $\rho(I_2)=\dfrac{6}{5}$, for all $r\ge 2$, $I_2^{(2r-1)} \subseteq I_2^r$, the above containment comes from the fact that $$\alpha(I_2^{(2r-1)}) =5r-2 \ge r-1+3r = r-1+\omega(I_2^r)$$ for all $r\ge 2$. The case $r=1$ is obvious.
\end{proof}

\begin{remark}
The above corollary gives a proof for the case $D_3$ in \cite[Proposition 6.3]{DrabkinSeceleanu}.
\end{remark}

We end this section by calculating $\beta(I_2^{(m)})$. 

\begin{proposition}
For all $m\ge 1$, $\beta(I_2^{(m)})=3m$ and $\omega(I_2^{(m)}) \ge 3m$.
\end{proposition}
\begin{proof}
The proof is the same as that of the case where $n\ge 3$. First, since $I_2^m \subseteq I_2^{(m)}$, we have that $\beta(I_2^{(m)}) \le 3m$ for all $m \ge 1$. On the other hand, recall that each line $L_j$ in the configuration contains exactly $3$ points of the configuration. Thus, for any $m\ge 1$ and for any $f \in [I_n^{(m)}]_t$ where $t<3m$, intersecting any line $L_j$ with the variety defined by $f$, by Bezout's theorem, since $\deg(f)\deg(L_j)<3m$, $L_j$ is a component of the variety of $f$. Therefore, for any $t<3m$, $L_j$ is a component of the zero locus of $[I_2^{(m)}]_t$, i.e, $\beta(I_2^{(m)}) \ge 3m$.
\end{proof}
 
\begin{remark}
As in the case where $n\ge 3$, Macaulay2 \cite{M2} suggests that $\omega(I_2^{(m)}) = 3m$. It is interesting to know if $\omega(J^{(m)})=\beta(J^{(m)})$ hold for what radical ideal of points $J$ in general. As suggested by the referee, the answer is no in general. Consider $8$ general points in $\PP ^2$ with defining ideal $I$. Then there are $2$ cubics among the generators of $I$, which form a regular sequence, i.e., intersect in $9$ points. By Cayley-Bacharach theorem, since any cubics containing $8$ points also contains the $9$-th point, one must use a form of degree at least $4$ (and in fact to be $4$) to exclude the $9$-th point from the defining ideal of $8$ points. Thus $\beta(I)=3$ but $\omega(I) \ge 4$.
\end{remark} 

\bibliographystyle{alpha}
\bibliography{References}  

\end{document}